\documentclass[11pt]{article}
\usepackage{latexsym,amsthm,verbatim,ifthen,graphicx,color,mathrsfs}
\usepackage[all]{xy}
\usepackage{color}

\usepackage{tikz-cd}
\usepackage{mathtools}

\usepackage{color}

\usepackage{lscape}
 \usepackage[ansinew]{inputenc}
 \usepackage{multicol}
 \usepackage[all]{xy}\xyoption{poly}
 \usepackage{mathrsfs}
\usepackage[psamsfonts]{eucal}
\usepackage{amssymb, amsmath}
\usepackage{geometry}
\usepackage{enumerate}
\usepackage{float}

\newtheorem{theorem}{Theorem}[section]

\newtheorem{corollary}[theorem]{Corollary}
 \newtheorem{lemma}[theorem]{Lemma}
 \newtheorem{proposition}[theorem]{Proposition}
 \theoremstyle{definition}
 \newtheorem{definition}[theorem]{Definition}

 \newtheorem{rem}[theorem]{Remark}
 \newtheorem{ex}[theorem]{Example}

\def\cG{{\mathscr G}}
\def\gauge{\,\cG\,}

\def\quil{{\mathscr L}}
\def\quilc{\widehat{\mathscr L}}
\def\quic{{\mathscr C}}

\def\esp#1{{\quad\text{#1}\quad}}

\newcommand{\bz}{\mathbb Z}

\def\bq{\mathbb{Q}}

\def\C{\mathbb{C}}

\def\Q{\mathbb{Q}}

\def\im{{\rm Im\,}}
\def\hL{{\widehat{\mathbb L}}}

\def\des{{s^{-1}}}
\newcommand{\fib}{\operatorname{{\cal F}\text{\rm ib}}}

\def\timest{\,\widetilde\times\,}

    \newcommand{\lasu}{{\mathfrak{L}}}
    
     \newcommand{\der}{{\rm Der}}
      \newcommand{\derr}{{\rm Der}}
       \newcommand{\derge}{{\mathcal{D}er}}
      \newcommand{\derpi}{{\rm Der}^{\Pi}}
      
       \newcommand{\derko}{{\rm Der}^{\mathcal K}_0}

 \newcommand{\lib }{\mathbb{L}}
  \newcommand{\calge }{{\cal G}}
  
    \newcommand{\calk }{{\cal K}}

 \newcommand{\fil}{\operatorname{{\mathcal F}}}
  \newcommand{\gil}{\operatorname{{\cal G}}}

  \newcommand{\Hom}{\operatorname{\text{\rm Hom}}}

\newcommand{\catdgl}{\operatorname{{\bf dgl}}}

\newcommand{\catcdgl}{\operatorname{{\bf cdgl}}}
\newcommand{\catcdgle}{\operatorname{{\bf cdgl^{\Delta}}}}

\newcommand{\catss}{\operatorname{{\bf sset}}}

\newcommand{\cale}{{\mathcal E}}

 \newcommand{\map}{\operatorname{{\rm map}}}
 \newcommand{\fmap}{\operatorname{{\mathfrak{map}}}}
  \newcommand{\Map}{\operatorname{{\rm map}}}
    \newcommand{\Der}{\operatorname{{\rm Der}}}
    \newcommand{\id}{\operatorname{{\rm id}}}
        \newcommand{\ad}{\operatorname{{\rm ad}}}

 \newcommand{\MC}{\operatorname{{\rm MC}}}
  \newcommand{\aut}{\operatorname{{\rm aut}}}

     \newcommand{\ev}{\operatorname{{\rm ev}}}
\newcommand{\mc}{{\MC}}

\newcommand{\catcdgc}{\operatorname{{\bf cdgc}}}

  \newcommand{\otimesc}{\widehat{\otimes}}

   \newcommand{\libc}{{\widehat\lib}}

\begin{document}

\title{Rational homotopy type of relative universal fibrations}

\author{Yves F\'elix, Mario Fuentes and Aniceto Murillo\footnote{The  authors have been partially supported by the  MICINN grant PID2020-118753GB-I00 of the Spanish Government and the Junta de Andaluc{\'\i}a grant ProyExcel-00827.}}

\maketitle

\begin{abstract} For any group $G$ of self homotopy equivalences of the finite nilpotent complex $X$, acting nilpotently on its homology,   and for any nilpotent subcomplex $A$, we prove that the universal fibration $$
X \longrightarrow B(*,\aut^{A}_G(X),X)\longrightarrow B\aut^{A}_G(X),
$$
which classifies $A$-fibrations  for which the image of the  $A$-holonomy action lies in $G$, has a Lie model of the form
$$
L\longrightarrow L\timest\derge^ML\longrightarrow\derge^ML
$$
in which: $M\hookrightarrow L$ is a Lie model of $A\hookrightarrow X$ and $\derge^ML$ is a connected complete differential graded Lie algebra of   derivations of  $L$ which vanish on $M$. The rational homotopy type of extended relative mapping fibrations is also similarly characterized.
\end{abstract}

\section*{Introduction}

For a given  finite pointed complex $X$, evaluating at the base point each self map of $X$ provides a fibration
\begin{equation}\label{cero}
\map^*(X,X)\longrightarrow \map(X,X)\stackrel{\ev}{\longrightarrow} X
\end{equation}
whose fibre at the base point consists of self pointed maps of $X$. Restricting this sequence to the components of homotopy equivalences produces a sequence of fibrations (each two maps constitute a fibration sequence)
\begin{equation}\label{primera}
\aut^*(X)\to\aut(X)\stackrel{\ev}{\to} X \to B\aut^*(X)\to  B\aut(X),
\end{equation}
where $\aut(X)$ denotes the topological monoid of self homotopy equivalences of $X$, $\aut^*(X)$ is the submonoid of pointed equivalences and $B$ stands for classifying space. The last fibration   classifies all fibrations sequences of fibre $X$ \cite{may,sta}.
Moreover, the pullback of this fibration over itself provides the universal fibration in the pointed setting,
\begin{equation}\label{segunda}
X\longrightarrow B(*,\aut^*(X), X)\longrightarrow B\aut^*(X),
\end{equation}
where $B(*,\aut^{*}(X),X)$ denotes the geometric bar construction \cite[Proposition 7.8]{may}. This sequence classifies pointed fibrations with fibre $X$ \cite{may}, i.e., fibrations of pointed spaces endowed with a pointed section.

The above can be stated more generally: given a subgroup $G\subset\cale(X)$ of homotopy classes of self homotopy equivalences, we will denote by $G^*\subset \cale^*(X)$ the subgroup of pointed classes of self homotopy equivalences whose free classes lie in $G$. Denote respectively  by $\aut^*_G(X)\subset \aut^*(X)$ and $\aut_G(X)\subset\aut(X)$ the submonoids of pointed and free equivalences whose free homotopy classes belong to $G$. That is $\pi_0\aut_G^*(X)=G^*$ and $\pi_0\aut_G(X)=G$. The {\em $G$-versions} of (\ref{primera}) and (\ref{segunda}) are, respectively, the fibration sequences
\begin{equation}\label{primeraprima}
\aut^*_G(X)\to\aut_G(X)\stackrel{\ev}{\to} X \to B\aut_G^*(X)\to  B\aut_G(X),
\end{equation}
and
\begin{equation}\label{segundaprima}
X\longrightarrow B(*,\aut^*_G(X), X)\longrightarrow B\aut^*_G(X).
\end{equation}
The last  sequence is also universal as  it classifies all pointed fibrations with fibre $X$ for which the image of the pointed holonomy action $\pi_1(X)\to \pi_0\aut^* (X)$ lies in $G^*$ \cite[Theorem 3.2]{fuen}.

In \cite[Theorem 0.2]{fefuenmu0} the authors build a sequence of complete differential graded Lie algebras (cdgl's from now on) which models this universal fibration in the modern sense of \cite{bufemutan0}:   whenever $X$ is nilpotent and $G$ acts nilpotently on the homology groups of $X$,  the geometrical realization of this cdgl sequence has the rational homotopy type of (\ref{segundaprima}). In particular, whenever  $X$ is simply connected, $G$ is the trivial group, and focusing only in the base of this sequence, we recover the classical result of Tanr\'e \cite[\S VII]{tan} modeling $B\aut^*_1(X)$ by means of derivations of the Quillen model of $X$.

In this paper we consider the relative setting: let $A$ be a subcomplex of the finite complex $X$ and consider the relative analogue of (\ref{segunda}) which is the fibration sequence
\begin{equation}\label{cuarta}
X \longrightarrow B(*,\aut^{A}(X),X)\longrightarrow B\aut^{A}(X),
\end{equation}
where $\aut^A(X)$ is the topological monoid of homotopy equivalences that fix $A$. This sequence classifies $A$-fibrations of fibre $X$ \cite[Theorem B.4]{heplah}, see \S1 for a precise definition and properties satisfied by these fibrations.

More generally, given a subgroup $G\subset\cale(X)$ denote by $G^A\subset \cale^A(X)$ the subgroup of relative homotopy classes of homotopy equivalences of $X$ which fix $A$ whose free classes live in $G$. In the topological side, denote by $\aut_G^A(X)\subset \aut^A(X)$ the submonoid of self homotopy equivalences whose relative homotopy classes lie in $G^A$. As a general picture, the group morphisms induced at the path components by the sequence
$$
\aut^A_G(X)\hookrightarrow \aut^*_G(X)\hookrightarrow \aut_G(X)
$$
is the subsequence
$$
G^A\longrightarrow G^*\longrightarrow G\quad\text{of}\quad
\cale^A(X)\longrightarrow \cale^*(X)\longrightarrow \cale(X).
$$
The $G$-version of (\ref{cuarta}) is the fibration sequence
$$
X \longrightarrow B(*,\aut^{A}_G(X),X)\longrightarrow B\aut^{A}_G(X),
$$
which classifies $A$-fibrations of fiber $X$ for which the image of the  $A$-holonomy action $\pi_1(B)\to\cale^A(X)$ lies in $G$.

Assume $A$ and $X$ nilpotent and let $G\subset\cale(X)$ be a subgroup acting nilpotently on the  homology groups of $X$. On the algebraic side let $j\colon M\hookrightarrow L$ be an inclusion of cdgl's modeling the inclusion $A\hookrightarrow X$ and denote by $\derr^ML\subset \derr L$ the sub dgl of derivations of $L$ which vanish at $M$ (see \S1 for a short but thorough compendium of the homotopy theory of cdgl's).  Although these dgl's are not complete in general, there is a Malcev $\bq$-complete subgroup $\aut_\calge^M(L)$ of the group $\aut(L)$ of automorphisms of $L$ whose logarithm   provides a complete Lie algebra $\derge_0^M L\subset \derr^M_0L$ of derivations of degree $0$. Extend this complete Lie algebra to a connected cdgl $\derge^ML$  by adding all positive derivations of $\derr^ML$. Then, see Theorem \ref{main1}:

\begin{theorem}\label{prin1} The twisted product
$$
L\longrightarrow L\timest\derge^ML\longrightarrow\derge^ML
$$
in which both terms are sub dgl's and $[\theta,x]=\theta(x)$ for any $\theta\in\derge^ML$ and any $x\in L$, is a cdgl fibration whose realization is  homotopy equivalent to the classifying fibration
$$
X_\bq \longrightarrow B(*,\aut^{A_\bq}_{G_\bq}(X_\bq),X_\bq)\longrightarrow B\aut^{A_\bq}_{G_\bq}(X_\bq).
$$
\end{theorem}

In the particular case in which $A$ and $X$ are simply connected, $G$ is the trivial group, and focusing in the base of this fibration, we retrieve the main result of \cite{bersa} exhibiting the simply connected cover of $\der^ML$ as a Lie model of $\aut_1^A(X)$ which is the universal cover of $\aut^A(X)$.

From Theorem \ref{prin1} we deduce, see Theorem \ref{teofin} and Corollary \ref{corofin}:

\begin{theorem}\label{prin2} The realization of the sequence of inclusions
$$
\derge^M L\hookrightarrow\derge L\hookrightarrow  \derge L\timest sL
$$
is homotopy equivalent to
$$
B\aut^{A_\bq}_{G_\bq}(X_\bq)\longrightarrow B\aut^*_{G_\bq}(X_\bq)\longrightarrow B\aut_{G_\bq}(X_\bq).
$$
\end{theorem}

It is important to remark that the proof of Theorem \ref{prin1} cannot rely on a relative extension of the approach used in \cite[Theorem 0.2]{fefuenmu0} where  the classifying property of (\ref{primeraprima}) and a previous description of its rational homotopy were needed. Indeed, the relative version of (\ref{primeraprima}) lack of classifying properties and the same strategy cannot be followed. Nevertheless, here we go the opposite way:

The analogue of (\ref{cero}) in the relative, free case, is the fibration
$$
\map^A(X,X)\longrightarrow\map(X,X)\stackrel{i^*}{\longrightarrow}\map(A,X)
$$
where $i\colon A\hookrightarrow X$ is the inclusion and $\map^A(X,X)$ denotes the fibre at $i$, i.e.,  the self maps of $X$ preserving $A$. The restriction of this sequence to the components  of  self homotopy equivalences of $X$ provides a sequence of fibrations
$$
\aut^A(X)\longrightarrow \aut(X)\stackrel{i^*}{\longrightarrow}\fmap(A,X)\longrightarrow B\aut^A(X)\longrightarrow B\aut(X),
$$
where $\fmap(A,X)\subset \map(A,X)$ consists of those  maps which can be extended to homotopy equivalences of $X$.

For a given subgroup $G\subset\cale(X)$, the $G$-version of this sequence reduces to
$$
\aut^A_{G}(X)\longrightarrow \aut_G(X)\stackrel{i^*}{\longrightarrow}\fmap_G(A,X)\longrightarrow B\aut^A_{G}(X)\longrightarrow B\aut_G(X),
$$
where $\fmap_G(A,X)\subset\fmap(A,X)$ is the subspace of those maps which can be extended to homotopy equivalences whose homotopy classes belong to $G$.

On the Lie side,   the inclusion $M\hookrightarrow L$   induces a short cdgl sequence of convolution Lie algebras (see \S2  for a brief compendium of these Lie algebras, and Remark \ref{coideal} for this precise statement)
$$
\Hom^M(\quic(L),L)\longrightarrow  \Hom(\quic(L),L)\longrightarrow\Hom(\quic(M),L),
$$
which is extended on the right by means of a natural dgl twisted product (see \S2) to obtain a quite simple dgl sequence formed by inclusions and projections,
$$
\begin{aligned}
\Hom^M&(\quic(L),L)\to  \Hom(\quic(L),L)\to\Hom(\quic(M),L)\to\\
&\to\Hom(\quic(M),L)\timest\derge L\timest sL\to \derge L\timest sL.
\end{aligned}
$$
Then, see Theorem \ref{final1} for a precise statement:
\begin{theorem}\label{prin3}
 The restriction to the appropriate components of this sequence produces a Lie model of the sequence
$$
\aut^{A_\bq}_{G_\bq}(X_\bq)\longrightarrow \aut_{G_\bq}(X_\bq)\stackrel{i_\bq^*}{\longrightarrow}\fmap_{G_\bq}(A_\bq,X_\bq)\longrightarrow B\aut^{A_\bq}_{G_\bq}(X_\bq)\longrightarrow \aut_{G_\bq}(X_\bq).
$$
\end{theorem}
Finally, in Theorem \ref{final2} we present an analogous result modeling the pointed version  of this sequence:
$$
\aut^{A_\bq}_{G_\bq}(X_\bq)\longrightarrow \aut_{G_\bq}^*(X_\bq)\stackrel{i_\bq^*}{\longrightarrow}\fmap_{G_\bq}^*(A_\bq,X_\bq)\longrightarrow B\aut^{A_\bq}_{G_\bq}(X_\bq)\longrightarrow \aut_{G_\bq}^*(X_\bq).
$$

To achieve our goal we first  present in \S2 some technical yet fundamental results relating convolution Lie algebras with Lie algebras of derivations.  It is also necessary, and we accomplish this in Section 3, to obtain explicit Lie models of relative mapping fibrations using convolution Lie algebras. With this groundwork laid, we proceed to prove   Theorems   \ref{prin1}, \ref{prin2} and \ref{prin3}, together with its pointed version,  in \S4.  Lastly, in Section 5 we give some examples/applications  of our main results. In particular, see Theorem \ref{teoejemplo}, we describe the group $\cale^{A_\bq}(X_\bq)$ whenever  $X$ is obtained by attaching a cell to the finite nilpotent CW-complex $A$.   It is however necessary to revisit the main aspects of the theory on which all of the above is based. This is done in \S1 where we include a subsection that outlines the counterpart of $A$-fibrations in the category $\catcdgl$ of complete differential graded Lie algebras.

\section{Preliminaries}
 We begin by setting some conventions of general nature: our notation will not distinguish a given category from the class of its objets. Every considered algebraic structure, say $L$, will be rational, $\bz$-graded and endowed with a (possibly null) differential, say $d$. If we want to specify such a differential  we write $(L,d)$. The {\em suspension} and {\em desuspension} of  such a structure  is denoted by $sL$ and $\des L$ respectively, i.e.,  $(sL)_n=L_{n-1}$ and $(\des L)_n=L_{n+1}$ for any $n\in\bz$.

Every considered topological space will be of the homotopy type of a CW-complex. Also, we often identify a simplicial set with the CW-complex obtained by its realization. Accordingly,  weak homotopy equivalences between simplicial sets are often refered to as homotopy equivalences. A subtle but important point here is that the realization of a simplicial mapping space is homotopy equivalent to the mapping space of the realization of the corresponding simplicial sets whenever these are finite complexes or their rationalizations.

\subsection{Homotopy theory of complete differential graded Lie algebras}

The general reference for most of what follows is the monograph \cite{bufemutan0}, see also \cite{bufemutan3,bufemutan1}. Here, we just emphasized the facts we need.

We denote by  $\catdgl$  the category of  differential graded Lie algebras, dgl's henceforward. A twisted product $L\timest L'$ of the dgl's $L$ and $L'$ consists of a dgl structure on $L\times L'$ for which $L\to L\timest L'\to L'$ is an exact dgl sequence.

Given a dgl $L$  we denote by  $\mc(L)$ the set of {\em Maurer-Cartan} (or just $\mc$) elements which are elements $a$ of degree $-1$  satisfying the Maurer-Cartan equation $da=-\frac{1}{2}[a,a]$. The perturbation of the differential $d$ of a given dgl $L$ by an MC-element $a$ is the differential $d_a= d+\ad_a$ where $\ad$ is the usual adjoint operator. The {\em component} of $L$ at $a$ is the connected (non negatively graded) sub dgl $L^a$ of $(L,d_a)$ given by
$$
L^a_p=\begin{cases} \ker d_a&\text{if $p=0$},\\ \,\,\,L_p&\text{if $p>0$}.\end{cases}
$$

A {\em complete differential graded Lie algebra}, cdgl henceforth,  is a dgl $L$ equipped with a  decreasing sequence of differential Lie ideals,
$$L=F^1\supset\dots \supset F^n\supset F^{n+1}\supset\dots$$ such that $[F^p,F^q]\subset F^{p+q}$ for $p,q\geq 1$, and  the  natural map
$$
L\stackrel{\cong}{\longrightarrow}\varprojlim_n L/F^n
$$
is a dgl isomorphism. Note that the lower central series of $L$,
$$
L^1\supset\dots\supset L^n\supset L^{n+1}\supset\dots,
$$
where $L^{1}= L$ and $L^{n}= [L, L^{n-1}]$ for $n>1$, is a filtration for any dgl which satisfies $L^n \subset F^n$ for any $n\ge 1$ and  any other filtration $\{F^n\}_{n\ge 1}$ of $L$. Unless specified otherwise, this is the considered filtration on a generic dgl. We denote by $\catcdgl$ the category of cdgl's whose morphisms are dgl morphisms which preserve the corresponding filtrations.

The  {\em completion} of a filtered dgl  is the dgl
$$
\widehat L=\varprojlim_nL/F^{n}
$$
which is complete with respect to the filtration
$
\widehat{F}^n=\ker ( \widehat L \to L/F^n)$.

Let $\lib(V)$ denote the free Lie algebra generated by the graded vector space $V$. Of particular importance in our theory, see \cite[\S3.2]{bufemutan0}, is the completion of a dgl of the form $(\lib(V),d)$ which is the cdgl
$$
\libc(V)=\varprojlim_n\lib(V)/\lib(V)^n.$$

 For a given cdgl $L$, the {\em gauge action} of  the group $L_0$, endowed  with the Baker-Campbell-Hausdorff product (BCH product henceforth), on the set
   $\mc(L)$ is defined by
$$
x\,\cG\, a=e^{\ad_x}(a)-\frac{e^{\ad_x}-1}{\ad_x}(d x)= \sum_{i\geq 0} \frac{\ad_x^i(a)}{i!} - \sum_{i\geq 0} \frac{\ad_x^i(dx)}{(i+1)!},\quad x\in L_0,\,a\in\mc(L).
$$
The orbit set $\mc(L)/\cG$  is denoted by $\widetilde\mc(L)$.

The homotopy theory of cdgl's is based in the pair of adjoint functors {\em (global) model} and {\em realization} \cite[Chapter 7]{bufemutan0},
$$
\xymatrix{ \catss& \catcdgl \ar@<1ex>[l]^(.50){\langle\,\cdot\,\rangle}
\ar@<1ex>[l];[]^(.50){\lasu}\\}.
$$
The category $\catcdgl$ inherits by transfer a cofibrantly generated model structure for which the above becomes  a Quillen pair \cite[Chapter 8]{bufemutan0}. To avoid confusion with Lie brackets, the set of homotopy classes of cdgl morphisms in the corresponding homotopy category  will be denoted by $\{\cdot\,,\cdot\}$.  In this new model structure fibrations are surjections in non negative degrees and weak equivalences between connected cdgl's are quasi-isomorphisms. Moreover, the induced functors in the respective homotopy categories extend the classical Quillen equivalences between rational homotopy types of simply connected simplicial sets and homotopy types of simply connected (i.e., positively graded) dgl's \cite{qui}. In this way, see \cite{fefuenmu2}, if $L$ is simply connected, $\langle L\rangle$ is weakly homotopy equivalent to the Quillen classical realization of $L$. In general, under no restriction, $\langle L\rangle$ is a strong deformation retract of $\mc_\bullet(L)$, the {\em Deligne-Getzler-Hinich groupoid} of $L$ (see \cite[\S11.4]{bufemutan0}). On the other hand,  If $X$ is a simply connected simplicial set of finite type and $a$ is any of its vertices, then \cite[Theorem 10.2]{bufemutan0}, $\lasu_X^a$ has the homotopy type of the Quillen model of $X$. Moreover, if $X$ is nilpotent, $\langle\lasu^a_X\rangle$ is  homotopy equivalent to $X_\bq$. In general, if $X$ is connected and of finite type, $\langle\lasu^a_X\rangle$ has the homotopy type of $\bq_\infty X$, the Bousfield-Kan $\bq$-completion of $X$ \cite{BK}.

For any cdgl $L$, the set $\mc(L)$ is precisely the set of $0$-simplices of $\langle L\rangle$ and, for any $a\in\mc(L)$, $\langle L\rangle^a$ has the homotopy type of $ \langle L^a\rangle$,  the path component of $\langle L\rangle$ containing  $a$. Moreover,
$$
      \langle L\rangle\simeq  \amalg_{a\in \widetilde{\mc}(L)} \,\, \langle L^a\rangle\quad\text{and}\quad \langle L\rangle\cong\langle (L,d_a)\rangle\quad\text{for any $a\in\mc(L)$}.
$$
Given a connected cdgl  $L$,  there are group isomorphisms
     $$
     \pi_n\langle L\rangle\cong H_{n-1}(L), \quad n\ge 1,
     $$
   where the group structure in $H_0(L)$ is considered with the BCH product.

     On the other hand, see  \cite[Chapter 7]{bufemutan0} for details, the global model $\lasu_X$ of a simplicial set $X$ completely reflects its simplicial structure. In particular,
 the  $0$-simplices of $X$ are the Maurer-Cartan elements of $\lasu_X$.

A {\em model} of a connected cdgl $L$ is a connected cdgl of the form $(\libc(V),d)$ together with a quasi-isomorphism
$$
(\libc(V),d)\stackrel{\simeq}{\longrightarrow} L.
$$
If $d$ is decomposable we say that $(\libc(V),d)$ is the {\em minimal model} of $L$ and is unique up to cdgl isomorphism.

\begin{definition}\label{minimamodel} Let $X$ be a connected simplicial set and let $a$ be any of its vertices. The {\em minimal model}  of $X$ is the minimal model $(\libc(V),d)$ of  $\lasu_X^a$.
\end{definition}

Then, $sV\cong \widetilde H_*(X;\bq)$ and, provided $X$ of finite type, $sH_*(\libc(V),d)\cong\pi_*(\bq_\infty X)$.

The {\em derivations}  of a given   dgl $L$ play an essential role in what follows. They form a dgl $\der L$ with the usual Lie bracket and differential $D=[d,-]$:
$$
[\theta,\eta]=\theta\circ\eta-(-1)^{|\theta||\eta|}\eta\circ\theta,
\quad
D\theta=d\circ\theta-(-1)^{|\theta|}\theta\circ d.
$$
However, see \cite[\S2]{fefuenmu1}, even for a generic simply connected dgl $L$, the connected cover of $\derr L$  may fail to be complete. Nevertheless, let $L=(\libc(V),d)$ be a free cdgl and choose a finite filtration of $V$ by graded subspaces,
$$V=V^0\supset V^1 \supset\dots \supset V^{q-1}\supset V^q=0.$$
Then  the connected dgl $\derge L\subset \der L$ given by
\begin{equation}\label{deresp}
\derge_k L=
 \left\{
 \begin{array}{cl}
 \der_k L, & \text{ if } k>0,
 \\
\theta\in \der_0 L, \textup{ such that } \theta(V^i)\subset V^{i+1}\oplus \hL^{\geq 2}(V), & \text{ if } k=0,
 \end{array}\right.
\end{equation}
is complete \cite[\S2]{fefuenmu1}.

On the other hand, given a dgl morphism $f\colon L\to L'$, denote by $\derr_f(L,L')$ the chain complex of {\em $f$-derivations} in which $\derr_f(L,L')_n$ are linear maps $\theta\colon L\to L'$ of degree $n$ such that
$$\theta [a,b] = [\theta (a), f(b)] + (-1)^{n  \vert a \vert} [f(a), \theta (b)],\quad a,b\in L.$$

\subsection{Relative fibrations}

 We say that a fibration $E\stackrel{p}{\to}B$ has fibre $F$ if this space is weakly homotopically equivalent to the fibre $p^{-1}(b)$ at the base point of $B$. More loosely, by a {\em fibration sequence} we mean a sequence of maps  $F\to E\to B$ whose composition is homotopic to a constant $b$ and the induced map from $F$ to the homotopy fiber of $E\to B$ over $b$ is a weak homotopy equivalence. The same  applies to the category $\catcdgl$.

 The geometrical realization $B(\cdot,\cdot,\cdot)$ of the simplicial topological space given by the geometric bar construction  plays an essential role in the classification of fibrations of different types and we refer to \cite{may} for its formal definition and main properties.

Recall from \cite[Appendix B]{heplah}, see also \cite[\S3.2]{fuen1}, that given a pair of CW complexes $(F,A)$, an $(F,A)$-fibration is a fibration of fibre $F$ over a connected space $B$ in which $A$ stays ``invariant'' along all the fibres. Precisely, it is a fibration $E\stackrel{p}{\to}B$ of fibre $F$ endowed with a map $\sigma\colon A\times B\to E$ over $B$ (i.e., $p\sigma=\text{proj}_B$)  such that, for each $b\in B$, the restriction of $\sigma$ to $A\times\{b\}$ is a homeomorphism onto its image and the following commutes
$$
\xymatrix{&A\ar@{^(->}[dl]\ar[dr]_\sigma&\\
F\ar[rr]^\simeq &&p^{-1}(b).\\}
$$
We say that $\sigma$ is an {\em $A$-section} of $p$.

As shown in \cite[Theorem B.4]{heplah}, $(F,A)$-fibrations are classified by means of the universal $(F,A)$-fibration
$$F\to B(\ast,\text{aut}^A(F),F) \to B\text{aut}^A(F).$$
More generally, see \cite[Theorem 2.39]{fuen1}, given any subgroup $G\subset \mathcal{E}(F)$, then
$$F\to B(\ast,\text{aut}^A_G(F),F)\to B\text{aut}^A_G(F)$$
is the universal fibration which classifies those $(F,A)$-fibrations such that the image of their \emph{relative} holonomy action $\pi_1(B)\to \cale^A(F)$  lies in  $G^A$.

The analogue of a relative fibration in the category $\catss$ of simplicial sets is the following: let $A\subset F$ be an inclusion of Kan complexes, and let $p\colon E\to B$ be a Kan fibration with fiber homotopy equivalent to $F$. We say that $p$ is a $(F,A)$-fibration if there exists a (degreewise) injective map
$\sigma\colon A\times B\to E$ over $B$, i.e.,
 $p \sigma=\text{proj}_B\colon A\times B\to B$.

As the geometric realization of simplicial sets takes Kan fibrations to (Serre) fibrations and preserves cofibrations, one easily checks that the  realization of an $(F,A)$-fibration $p\colon E\to B$ in $\textbf{sset}$ is a $(|F|, |A|)$-fibration.

In particular, by the classification theorem of $(|F|, |A|)$-fibrations, $|p|$ is obtained as the  pullback of the universal  $(|F|, |A|)$-fibration along a map
$|B| \to B\text{aut}^{|A|}{|F|}$
which, by adjointness of the geometrical realization and the singular functors, corresponds to a simplicial map $B\to \text{Sing}(B\text{aut}^{|A|}{|F|})$. As there is no possible confusion, we denote   $\text{Sing}(B\text{aut}^{|A|}{|F|})$ simply by $B\text{aut}^A(F)$, and deduce that the Kan fibration
$$F\to B(\ast,\text{aut}^A(F),F) \to B\text{aut}^A(F)$$
is universal with respect to $(F,A)$-fibration sequences in $\catss$. More generally, an analogous argument produces, for any given subgroup $G\in\cale(F)$, a universal Kan fibration of simplicial sets
$$F\to B(\ast,\text{aut}^A_G(F),F)\to B\text{aut}^A_G(F)$$
classifying $(F,A)$-fibrations of simplicial sets for which the images of their relative holonomy actions lie in $G^A$.

Finally, we introduce the corresponding term in $\textbf{cdgl}$.

\begin{definition}\label{relfib} Let
$$L\longrightarrow L'\stackrel{p}{\longrightarrow} L''$$
be a short exact sequence of connected cdgl's and $M\subset L$ a subcdgl. We say that $p$ is an {\em $(L,M)$-fibration} if there is an injective morphism   $\sigma\colon  M\times L''\to L'$ such that $p\sigma\colon M\times L'\to L'$ is the projection. We call $\sigma$ an {\em $M$-section}.
\end{definition}

As the realization functor preserves fibrations (it is a Quillen right adjoint) and injective maps (by definition), it follows that the realization of the cdgl fibration $p$,
$$\langle L\rangle \longrightarrow \langle L'\rangle\stackrel{\langle p\rangle}{\longrightarrow} \langle L''\rangle,$$
is a $(\langle L\rangle,\langle M\rangle)$-fibration of simplicial sets.

\section{Convolution Lie algebras and Lie algebras of  relative\hfill\break derivations}

We denote by  $\catcdgc$  the category of cocommutative differential graded coalgebras (cdgc's from now on) which are  assumed to have a counit $\varepsilon\colon C\to\bq$ and a coaugmentation $\eta\colon\bq\to C$, so that $C\cong\overline C\oplus \bq$ for  $\overline{C}= \ker\varepsilon$, and the reduced diagonal  $\overline{\Delta}\colon \overline C\to \overline C\otimes\overline C$ is given  by
$
\overline{\Delta}x= \Delta x- (u\otimes x + x\otimes u)
$,
with $u=\eta(1)$.

Recall that the categories $\catdgl$ and $\catcdgc$ are related by the classical Quillen functors
$$
\xymatrix{
{\rm \bf cdgc} \ar@<0.75ex>[r]^-{\mathscr L} &\mathbf{dgl} \ar@<0.75ex>[l]^(0.49){{\mathscr C} }
}
$$
where
${\mathscr L}(C)=(\mathbb L(s^{-1}\overline{C}),d)$,  and  $d=d_1+d_2$ with
$d_1(s^{-1}c) = -s^{-1}dc$ and
$d_2(s^{-1}c) = \frac{1}{2} \sum_i (-1)^{\vert a_i\vert} [s^{-1}a_i, s^{-1}b_i]$, with $\overline{\Delta}c = \sum_i a_i\otimes b_i$.

On the other hand,
${\mathscr C}(L)=(\land (sL),d)$ is the free cocommutative coalgebra cogenerated by  $sL$, and
$d=d_1+d_2$ where, for instance,
$d_1(sv) = -sdv$ and $d_2(sv\land sw)= (-1)^{\vert sv\vert} s[v,w]$.

The {\em convolution Lie algebra} associated to  any cdgc $C$ and any dgl $L$ is the dgl of linear maps $\Hom(C,L)$ with the usual differential $Df=d\circ f-(-1)^{|f|}f\circ d$, and the  Lie bracket,
$$
[f,g]=[\,\,,\,]\circ (f\otimes g)\circ\Delta.
$$
If we replace the diagonal by the reduced one we obtain  the {\em reduced} convolution Lie algebra $\Hom(\overline C,L)$.

 If $L$ is complete, then $\Hom(C,L)$ and $\Hom(\overline C,L)$ are also complete. Moreover, their homotopy type are invariants of the homotopy type of $C$ (see \cite[Remark 4.1]{fefuenmu0}).

There is also a natural twisted product, heavily used in what follows,
\begin{equation}\label{mainrelation}
\Hom(C,L)\timest \derr L
\end{equation}
where both terms are sub dgl's and
$
[\theta,f]=\theta f$ for any $\theta\in\derr L$ and any $ f\in \Hom(C,L)$.

\begin{rem}\label{coideal}  Note that If $B\subset C$ is a cdgc coideal there is a natural short exact sequence of convolution dgl's
$$
0\to \Hom(C/B,L)\longrightarrow\Hom(C,L)\longrightarrow \Hom(B,L)\to 0.
$$
This sequence does not split in general unless we choose $B=\bq$  in which case we have a dgl isomorphism,
\begin{equation}\label{iso}
\Hom(C,L)\cong L\timest \Hom(\overline C,L),
\end{equation}
where $L$ is identified with $\Hom(\bq,L)$, both factors of the twisted product are sub dgl's, and $[x,f]=\ad_x\circ f$ for  $x\in L$ and $f\in \Hom(\overline C,L)$. Under this splitting one exhibits
$$
\Hom(\overline C,L)\timest \derr L
$$
as sub dgl of (\ref{mainrelation}).

 In the reduced case there is also a short exact sequence of dgl's
\begin{equation}\label{finali}
0\to \Hom(\overline C/\overline B,L)\longrightarrow\Hom(\overline C,L)\longrightarrow \Hom(\overline B,L)\to 0.
\end{equation}
Note that
$$ \Hom(\overline C/\overline B,L)\cong \Hom( C/ B,L)
\cong \Hom^B(C,L)
$$
where the latter  denotes the sub dgl of $\Hom(C,L)$  consisting on maps that are trivial on $B$.
\end{rem}

\begin{rem}\label{impo2}
It is well known that
$$
\mc\bigl(\Hom(C,L)\bigr)\cong \Hom_{\catcdgc}\bigl(C,\quic(L)\bigr).
$$
Indeed,  consider   the degree $-1$ linear map,
\begin{equation}\label{q}
q\colon \quic(L)\to L\quad q(1)=q(\land^{\ge 2}sL)=0,\quad q(sx)=x,\quad x\in L.
\end{equation}
and denote in the same way its restriction $q\colon\overline\quic(L)\to L$.
It is then straightforward to check that the map
$$
\Hom_{\catcdgc}\bigl(C,\quic(L)\bigr)\to \mc\bigl(\Hom(C,L)\bigr),\quad \psi\mapsto q\psi,
$$
is a bijection. In the same way,
$$
\mc\bigl(\Hom(\overline C,L)\bigr)\cong \Hom_{\catcdgc}\bigl(\overline C,\overline \quic(L)\bigr).
$$
In particular, whenever $C=\quic(L')$ for some  dgl $L'$, any dgl morphism $\varphi\colon L'\to L$ produces the MC elements
$$
\varphi q=q\quic(\varphi)\in \mc\bigl(\Hom(\quic(L'),L)\bigr)\quad\text{and}\quad  \varphi q=q\overline\quic(\varphi)\in \mc\bigl(\Hom(\overline\quic(L'),L)\bigr).
$$
\end{rem}

In this section we describe the close correlation  existing between convolution Lie algebras and Lie algebras of derivations. The geometrical meaning of this association will be revealed in \S4.

 \smallskip

For any sub dgl $M\subset L$ we  consider the sub dgl's  of $\Hom(\quic(L),L)\timest \derr L$ given by
$$
\Hom(\overline\quic(L),L)\timest \derr^ML\quad\text{and}\quad
 \Hom(\quic(L),L)\timest \derr^M L.
$$
To avoid an excessive notation from now on, we abuse it and write $\Hom^M(\quic(L),L')$ instead of  $\Hom^{\quic(M)}(\quic(L),L')$.

Assume here after that $M\subset L$ is an inclusion of connected free cdgl's of the form
$$
M=(\libc(U),d)\stackrel{j}{\hookrightarrow}(\libc(U\oplus V),d)=L,
$$
and let
$$
\begin{aligned}
&(\Hom(\quic(L),L)\timest \derr^M L,D_q)\stackrel{\rho}{\longrightarrow} (\Hom(\quic(M),L),D_{jq}),\\
&(\Hom(\overline\quic(L),L)\timest \derr^M L,D_q)\stackrel{\overline\rho}{\longrightarrow} (\Hom(\overline\quic(M),L),D_{jq}),
\end{aligned}
$$
be the perturbed  surjections which vanish on the second factor and is the projection on the first according to Remark \ref{coideal}. Then:

\begin{proposition} \label{fundamental} Both $\rho$ and $\overline \rho$ are dgl quasi-isomorphisms.
\end{proposition}

\begin{lemma}\label{lemader} Let $\varphi\colon L'\to L$ be a dgl morphism in which $L'$ is free complete. Then, the chain map
$$
\Phi\colon s^{-1}\der_\varphi(L',L)\longrightarrow (\Hom(\overline\quic(L'),L), D_{\varphi q}),\quad \Phi(s^{-1}\theta)=\theta q,
$$
is a quasi-isomorphism.
\end{lemma}
\begin{proof} On the one hand the map induced by the unit $\alpha\colon \quil\quic(L')\stackrel{\simeq}{\to} L'$,
$$
s^{-1}\der_\varphi(L',L)\stackrel{\simeq}{\longrightarrow}s^{-1}\der_{\varphi\alpha}(\quil\quic(L'),L)
$$
is a quasi-isomorphism of chain complexes by virtue of \cite[Lemma 6]{bufemu0}.

On the other hand consider the isomorphism of chain complexes, explicitly detailed for instance in the proof of \cite[Theorem 5.1]{fefuenmu0},
$$
 s^{-1}\der_{\varphi\alpha}(\quil\quic(L'),L)\stackrel{\cong}{\longrightarrow} (\Hom(\overline\quic(L'),L), D_{\varphi q}),\quad(s^{-1}\eta)(c)=(-1)^{\eta}\eta(s^{-1} c).
$$
The composition of these maps is precisely $\Phi$.
\end{proof}

\begin{rem}\label{sucesionexacta} As a result we have the following  quasi-isomorphism of short exact sequences
$$
\xymatrix{
(\Hom^M(\quic(L),L),D_{q})\ar[r]&
(\Hom(\overline\quic(L),L),D_{q})\ar[r]&(\Hom(\overline\quic(M),L),D_{jq})\\
s^{-1}\derr^ML\ar[u]^\simeq_\Phi\ar[r]&
s^{-1}\derr L\ar[u]^\simeq_\Phi\ar[r]&
s^{-1}\derr_j(M,L)\ar[u]^\simeq_{\Phi}\\}
$$
where the upper row is (\ref{finali}) of Remark \ref{coideal} perturbed by $q$ and $jq$ respectively.
\end{rem}

\begin{proof}[Proof of Proposition \ref{fundamental}]
One easily checks that both maps commute with bracket and differentials. We first deal with the reduced case and consider the commutative diagram
$$
\xymatrix{
(\Hom(\overline\quic(L),L)\timest \derr^M (L),D_q)\ar[r]^(.57){\overline\rho}&(\Hom(\overline\quic(M),L),D_{jq})\\
(s^{-1}\derr(L)\timest \derr^ML,\widetilde D)\ar[r]^(.57)r\ar[u]^\simeq_{\Phi\times \id}&
s^{-1}\derr_j(M,L)\ar[u]^\simeq_{\Phi}\\}
$$
where: $r$ vanish on $\derr^ML$ and is the projection on the first factor according to Remark \ref{sucesionexacta};  $\widetilde D$ restricts to $s^{-1}D$ in $s^{-1}\derr L$; and finally, $\widetilde D\theta=D\theta+s^{-1}\theta$ for $\theta\in \derr^ML$. Then  $\ker r=(s^{-1}\derr^ML\timest \derr^ML,\widetilde D)$ where again $\widetilde D$ restricts to $s^{-1}D$ in $s^{-1}\derr L$ and $\widetilde D\theta=D\theta+s^{-1}\theta$ for $\theta\in \derr^ML$. This is trivially acyclic so $r$, and hence $\rho$, are quasi-isomorphisms.

For the unreduced case, simply note that, by the isomorphism (\ref{iso}) of Remark \ref{coideal}, the underline chain map of $\rho$ is identified with,
$$
\id\times \overline\rho \colon (L\times \Hom(\overline\quic(L),L)\times \der^ML, D_q)\longrightarrow (L\times\Hom(\overline\quic(M),L), D_{jq})
$$
where $D_qx=dx-(-1)^{|x|}\ad_x q$ and $D_{jq}x=dx-(-1)^{|x|}\ad_x j q$ for $x\in L$.
\end{proof}
Next, denote by $\Hom^{\overline\quic(M)}(\quic(L),L)$ the subdgl of $\Hom(\quic(L),L)$ of maps which vanish on $\overline\quic(M)$. Note that
\begin{equation}\label{split}
\Hom^{\overline\quic(M)}(\quic(L),L)\cong L\timest \Hom^M(\quic(L),L)
\end{equation}
and consider then the subdgl of $\Hom(\quic(L),L)\timest \derr L$  given by
$$
\Hom^{\overline\quic(M)}(\quic(L),L)\timest \der^ML\cong L\timest \Hom^M(\quic(L),L)\timest \der^ML.
$$
Even though  this subdgl is not closed for the perturbed $D_q$ we may endow it with the differential $\widetilde D$ which restricts to $D_q$ in $\Hom^M(\quic(L),L)\timest\der^ML$ and
$$\widetilde Dx=dx-(-1)^{|x|}\ad_x^Mq,\quad x\in L.
$$
Here $\ad_x^M\in\der^ML$ denotes the derivation which vanishes on $M$ and is $\ad_x$ on $V$. Then:
\begin{proposition} The projection
$$
\eta\colon(L\timest \Hom^M(\quic(L),L)\timest \der^ML,\widetilde D)\stackrel{\simeq}{\longrightarrow} L
$$
is a dgl quasi-isomorphism.
\end{proposition}\label{preseccion}
\begin{proof}A mere computation checks that $\eta$ a dgl morphism. Consider the commutative diagram of chain complexes
$$
\xymatrix{(L\times \Hom^M(\quic(L),L) \times  \der^ML,\widetilde D)\ar[r]^(.80){\eta}&L\\
(L\times s^{-1}\der^ML \times \der^ML,\widetilde D)\ar[ur]^{\eta'}\ar[u]^{\id\times\Phi\times\id}_{\simeq}&\\}
$$
 in which $\eta'$ is also the projection and where,
 in the bottom complex: $\widetilde D$ restricts to $s^{-1}D$ in $s^{-1}\der^ML$, $\widetilde Dx=dx+s^{-1}\ad_x^M$ for  $x\in L$ and $\widetilde D\theta=D\theta+s^{-1}\theta$ for $\theta\in \der^ML$. The kernel of $\eta'$ is trivially acyclic so this map, and thus $\eta$, are quasi-isomorphisms.
\end{proof}
A short computation proves the following:
 \begin{corollary}\label{seccion}
 The map
 $$
 \sigma\colon L\stackrel{\simeq}{\longrightarrow} (L\times s^{-1}\der^ML \times \der^ML,\widetilde D),\quad \sigma(x)=x-\ad_x^M,
 $$
 is a dgl  section of $\eta$ and thus, it is a quasi-isomorphism.
 \hfill$\square$
 \end{corollary}

\begin{proposition}\label{varsigmareducida}
The injection on the second factor
$$
\overline\kappa\colon \der^ML\stackrel{\simeq}{\hookrightarrow }(\Hom(\overline\quic(M),L)\timest\der L,D_{jq})
$$
is a dgl quasi-isomorphism.
\end{proposition}

\begin{proof} A simple computation shows that $\overline\kappa$ is a dgl morphism. Consider the commutative diagram of chain complexes

$$
\xymatrix{&(s^{-1}\der_j(M,L)\times\der L,\widetilde D)\ar[d]^{\Phi\times\id}_{\simeq}\\
\der^ML\ar@{^(->}[ur]^{\overline\kappa'}\ar@{^(->}[r]_(.27){\overline\kappa}&(\Hom(\overline\quic(M),L)\timest\der L,D_{jq})\\}
$$
where $\overline\kappa'$ is  the inclusion on the second factor, $\widetilde D$ restricts to $s^{-1}D$ on $s^{-1}\derr_j(M,L)$ and $\widetilde D\theta=D \theta+s^{-1}(\theta j)$ for $\theta\in\der L$. Observe that
$$
\begin{aligned}
(s^{-1}\der_j(M,L)\times\der L)/\der^ML&\cong s^{-1}\der_j(M,L)\times(\der L/\der^ML)\\
&\cong s^{-1}\der_j(M,L)\times \der_j(M,L).
\end{aligned}
$$
Moreover, the differential $\widetilde D$ inherited in the last chain complex is $s^{-1}D$ on $s^{-1}\der_j(M,L)$ and $\widetilde D \theta=D\theta +s^{-1}\theta$ for $\theta\in\der_j(M,L)$.
This is acyclic so $\overline\kappa'$, and thus $\overline\kappa$, are quasi-isomorphisms.
\end{proof}

For the unreduced version of this result extend the dgl
$
\Hom(\quic(M),L)\timest\der L
$
 to the twisted dgl
$$
\Hom(\quic(M),L)\timest\der L\timest sL
$$
where $sL$ is an abelian sub Lie algebra and
$$
Dsx=-sdx+ad_x-x,\quad[\theta,sx]=(-1)^{|\theta|}s\theta(x),\quad[f,sx]=0,
$$
with $sx\in sL$, $\theta\in\der L$ and $f\in\Hom(\quic(M),L)$. Here, $x\in\Hom(\quic(M),L)$ denotes the map which is zero on $\overline\quic(M)$ and takes 1 to $x$.

\begin{proposition}\label{proponueva} The inclusion on the second factor
$$
\kappa\colon\der^ML\stackrel{\simeq}{\hookrightarrow }(\Hom(\quic(M),L)\timest\der L\timest sL,D_{jq}).
$$
is a dgl quasi-isomorphism.
\end{proposition}
\begin{proof}
Under the splitting  in (\ref{iso}) this map is,
$$
\kappa\colon\der^ML\stackrel{\simeq}{\hookrightarrow }(\Hom(\overline\quic(M),L)\timest L\timest\der L\timest sL,D_{jq}).
$$
which is simply $\overline\kappa$ in which the codomain has been extended.

On the other hand observe that,
  as chain complexes,
$$
(\Hom(\overline\quic(M),L)\timest L\timest\der L\timest sL,D_{jq})/(\Hom(\overline\quic(M),L)\timest\der L,D_{jq})\cong (L\times sL,\widetilde d)
$$
where $\widetilde d=d$ on $L$ and $\widetilde d(sx)=-sdx-x$. This is then an acyclic complex and the result follows from Proposition \ref{varsigmareducida}.
\end{proof}

\section{Lie models of relative mapping fibrations}\label{relativogeneral}

For the rest of the section we fix a pointed inclusion $i\colon A{\hookrightarrow}X$  of connected, finite type, nilpotent complexes modeled by the inclusion $ M\hookrightarrow L$, and denote by $Y=\langle L'\rangle$ the realization of the connected cdgl $L'$. Note that $Y$ is canonically pointed as it has just one $0$-simplex. It is convenient to keep in mind that, whenever $L'$ is a Lie model of a connected complex $Z$ of finite type, then $Y\simeq Z^{\wedge}_\bq$.

\begin{proposition}\label{modelohomrelnil}
The realization of the maps induced by $ M\hookrightarrow L$,
$$
 \Hom(\quic(L),L')\to\Hom(\quic(M),L')\quad \text{and}\quad \Hom(\overline\quic(L),L')\to
 \Hom(\overline\quic(M),L'),
$$
are weakly homotopy equivalent to the relative mapping fibrations
$$i^*\colon\map(X,{Y})\to \map(A,{Y})\quad \text{and}\quad i^*\colon\map^*(X,{Y})\to \map^*(A,{Y}).
$$
\end{proposition}

\begin{proof} In the unpointed case, if $B$ is any cdga model of finite type of $X$,  \cite[Theorem 10.2]{bufemutan0} asserts that $L\simeq \quil(B^\sharp)$ and thus $\quic(L)\simeq\quic\quil(B^\sharp)\simeq B^\sharp$. Thus,
$$
B\otimesc L'\cong \Hom(B^\sharp,L')\cong\Hom(\quic(L),L').
$$
Hence,  by \cite[Theorem 12.18]{bufemutan0}, the realization of  $\Hom(\quic(L),L')$ is weakly equivalent to $\map(X,{Y})$. The naturality of this equivalence finishes the proof.
For the pointed case proceed analogously invoking this time \cite[Proposition 12.25]{bufemutan0}.
\end{proof}

We also provide a natural  Lie model of any non empty fibre  of the  free and pointed mapping fibrations at a given  map  $h\colon A\to{Y}$. Consider  the fibration sequences
\begin{equation}\label{fibrageneral}
F_{h}\to\map(X,{Y})\stackrel{i^*}{\to} \map(A,{Y}) \quad\text{and}\quad F_{h}\to\map^*(X,{Y})\stackrel{i^*}{\to} \map^*(A,{Y})
\end{equation}
where $F_{h}\subset\map(X,{Y})$ is the  non path-connected subspace of maps extending $h$, which is assume to be pointed in the based case.
As we assume this is not the empty space, there is a   cdgl morphism $\varphi\colon L\to L'$ whose restriction to $M$ models $h$. Following Remark \ref{coideal},
 $\Hom^M(\quic(L),L')\subset\Hom(\quic(L), L')$ denotes the sub cdgl of linear maps which vanish at $\quic(M)$.

\begin{proposition}\label{modelofibrageneral} The realization of $(\Hom^M(\quic(L),L'),D_{\varphi q})$  has the weak homotopy type of $F_h$.
\end{proposition}

\begin{proof}
 By Proposition \ref{modelohomrelnil}, as the realization of a cdgl does not depend on  perturbing differentials, the cdgl morphism
$$
(\Hom(\quic(L),L'),D_{\varphi q})\longrightarrow(\Hom(\quic(M),L'),D_{\varphi q})
$$
is a Lie model of $i^*$ where now, the trivial MC element in $(\Hom(\quic(M),L'),D_{\varphi q})$ represents $h$ (see Remark \ref{igual} below). Being surjective, this morphism is a cdgl fibration and thus, its kernel, which is precisely  $(\Hom^M(\quic(L),L'),D_{\varphi q})$, models $F_h$.
\end{proof}

\begin{corollary}\label{coroim}
The realization of the short exact cdgl sequences,
$$
(\Hom^M(\quic(L),L'),D_{\varphi q})\to(\Hom(\quic(L),L'),D_{\varphi q})\to(\Hom(\quic(M),L'),D_{\varphi q})
$$
and
$$
(\Hom^M(\quic(L),L'),D_{\varphi q})\to (\Hom(\overline\quic(L),L'),D_{\varphi q})\to(\Hom(\overline\quic(M),L'),D_{\varphi q})
$$
are weakly homotopy equivalent to
$$
F_{h}\to\map(X,{Y})\stackrel{i^*}{\to} \map(A,{Y}) \quad \text{and}\quad F_{h}\to\map^*(X,{Y})\stackrel{i^*}{\to} \map^*(A,{Y}).
$$
\hfill$\square$
\end{corollary}

The following, which in particular covers  the main results in \cite[\S2]{fefuenmu0} together with  Remark 4.5 in op. cit., provides essential observations:

\begin{rem}\label{igual} Note that
$$
\pi_0\Map(X, {Y})\cong \{L,L'\}/H_0(L')\cong \widetilde\mc\bigl(\Hom(\quic(L),L')\bigr)
$$
and
$$
\pi_0\Map^*(X, {Y})\cong\{L,L'\}\cong \widetilde\mc\bigl(\Hom(\overline\quic(L),L')\bigr).
$$
The same applies replacing $X$ and $L$ by $A$ and $M$.
We explicitly recall here these bijections:

On the one hand \cite[Corollary 2.2]{fefuenmu0}, the set $[X,Y]^*$ of pointed homotopy classes is in bijective correspondence with $\{L,L'\}$. Moreover, the homotopy class of a cdgl morphism $\varphi\colon L\to L'$, which corresponds to the pointed homotopy class of a pointed map $f\colon X\to Y$, is in turn identified to the gauge class of the MC element
$$
\varphi  q \colon\overline\quic(L)\longrightarrow L',
$$
with $q$ as  in (\ref{q}).

On the other side \cite[Corollary 2.5]{fefuenmu0}, the set $[X,Y]$ of free homotopy classes is identified with the orbit set $\{L,L'\}/H_0(L')$ where the action of $H_0(L')$ (with the BCH product) on $\{L,L'\}$ is   induced by the exponential: $a\cdot \varphi=\varphi\circ e^{\ad_a}$ for $a\in L_0$ and $\varphi\colon L\to L'$. Moreover, the class in this quotient of the cdgl morphism $\varphi\colon L\to L'$, which corresponds to the unbased homotopy class of $f$,  is identified to the gauge class of the  MC element
$$
\varphi q \colon \quic(L)\longrightarrow L'.
$$
\end{rem}

\begin{definition} Given a cdgl morphism $\gamma\colon M\to L'$ we denote by $\{L,L'\}^M_\gamma$ the set of homotopy classes relative to $M$ of cdgl morphisms which restrict to $\gamma$ in $M$. In particular, if $j\colon M\hookrightarrow L$ is the inclusion, $\{L,L\}^M$ will denote $\{L,L\}^M_j$.
\end{definition}

\begin{rem}\label{remim} Note that $\widetilde\mc\bigl(\Hom^M(\quic(L),L'),D_{\varphi q})\bigr)\cong\{L,L'\}^M_\varphi$. Indeed,
recall that, for any dgl $K$ and any $z\in\mc(K)$, the map $a\mapsto z-a$ produces a bijection $\mc(K)\cong\mc(K^z)$. In our setting, see Remark \ref{impo2}, this translates to
$$
\mc\bigl(\Hom(\quic(L),L'),D_{\varphi q})\bigr)=\{(\gamma-\varphi)q\,\,\text{with}\,\,\gamma\colon L\to L'\,\,\text{cdgl morphism}\}.
$$
In particular
$$
\mc\bigl(\Hom^M(\quic(L),L'),D_{\varphi q})\bigr)=\{(\gamma-\varphi)q\,\,\text{with}\,\,\gamma\colon L\to L'\,\,\text{such that $\gamma_{|_M}=\varphi_{|_M}$}\}.
$$
In other words, the $\mc$ elements of $\bigl(\Hom^M(\quic(L),L'),D_{\varphi q})\bigr)$ are identified with the set of cdgl morphisms $\gamma\colon L\to L'$ which restrict to $\varphi$ in $M$. Through this identification, and by Remark \ref{igual}, $\widetilde\mc\bigl(\Hom^M(\quic(L),L'),D_{\varphi q})\bigr)$ is in bijective correspondence with $\{L,L'\}^{M}_\varphi$.

Thus the different components of the fibre $F_h$  have the homotopy type of the realization of the components
$$
\Hom^M(\quic(L),L')^{\gamma q}=\Hom^M(\quic(L),L'),D_{\varphi q})^{(\gamma-\varphi)q}
$$
as $[\gamma]$ varies in $\{L,L'\}^{M}_\varphi$.
\end{rem}

\section{Lie models of relative universal fibrations}

Following the notation in \S3  we consider a pointed inclusion of nilpotent finite complexes $i\colon A{\hookrightarrow} X$. We fix $j\colon M\hookrightarrow L$ a Lie model of $i$ of the form
$
 (\libc(U),d)\hookrightarrow(\libc(U\oplus V),d)$ with $(\libc(U),d)$ and $(\libc(V),\overline d)$ minimal.

We consider $\aut(L)$ the group of automorphisms of $L$ and denote by    $\cale^*(L)\subset\{L,L\}$ the group $\aut(L)/\sim$  of homotopy classes of homotopy automorphisms of $L$. In the same way, if $\aut^M(L)$ is the subgroup of relative automorphisms of $\aut(L)$ which restricts to the inclusion on $M$, we denote by $\cale^M(L)\subset\{L,L\}^M$ the group  $\aut^M(L)/\sim_M$  of relative homotopy classes of relative automorphisms of $L$. Finally we also define $\cale(L)=\cale^*(L)/H_0(L)$ so that the natural sequence of group morphisms
$$
\cale^M(L)\longrightarrow \cale^*(L)\longrightarrow \cale(L)
$$
is identified to
$$
\cale^{A_\bq}(X_\bq)\longrightarrow \cale^*(X_\bq)\longrightarrow \cale(X_\bq).
$$
With the notation in \S3 we begin with:

\begin{theorem}\label{autofibra}
The fibrations
$$
\aut^{A_\bq}(X_\bq)\to \aut(X_\bq)\stackrel{i^*}{\to} \fmap(X_\bq,X_\bq)\,\,\text{and}\,\,\aut^{A_\bq}(X_\bq)\to \aut^*(X_\bq)\stackrel{i^*}{\to} \fmap^*(X_\bq,X_\bq)
$$
have, respectively, the homotopy type of
$$
\underset{[\gamma]\in\cale^M(L)}{\amalg} \langle \Hom^M(\quic(L),L)^{\gamma q}\rangle \to\underset{[\varphi]\in\cale(L)}{\amalg}\langle \Hom(\quic(L),L)^{\varphi q}\rangle \to\underset{[\varphi]\in\cale(L)}{\amalg}\langle \Hom(\quic(M),L)^{\varphi q}\rangle.
$$
and
$$
\underset{[\gamma]\in\cale^M(L)}{\amalg} \langle \Hom^M(\quic(L),L)^{\gamma q}\rangle \to\underset{[\varphi]\in\cale^*(L)}{\amalg}\langle \Hom(\overline\quic(L),L)^{\varphi q}\rangle \to\underset{[\varphi]\in\cale^*(L)}{\amalg}\langle \Hom(\overline\quic(M),L)^{\varphi q}\rangle.
$$
\end{theorem}
\begin{proof}
The fibrations in the statement are the restrictions of the unbased and pointed mapping fibrations in (\ref{fibrageneral}) to $\aut(X_\bq)\subset\map(X_\bq,X_\bq)$ and $\aut^*(X_\bq)\subset\map^*(X_\bq,X_\bq)$, where now the fibers are taking on $h=i_\bq\colon A_\bq\hookrightarrow X_\bq$.  The result follows then by restricting the realization of the Lie models of this fibrations in Corollary \ref{coroim} to the appropriate components. For it, take into account Remarks \ref{igual} and \ref{remim}.
\end{proof}
Fix now a subgroup $G\subset\cale(X)$ acting nilpotently on $H_*(X)$.  As thoroughly explained in \cite[\S7]{fefuenmu0}, its rationalization $G_\bq$ can be considered as a subgroup of $\cale(X_\bq)$ and via this identification the map $G\to G_\bq$, $[f]\mapsto[f_\bq]$, is also the rationalization.
\begin{definition}\label{explog} Let
$
\calge\subset \cale^*(L)$
be the subgroup of homotopy classes of automorphism of $L$ such that
$$
\calge/H_0(L)\cong G_\bq.
$$
Define accordingly the subgroup of $\aut(L)$ given by
 $$
 \aut_{\calge}(L)=\{\varphi\in\aut(L),\,[\varphi]\in \calge\}.
 $$
 In the relative setting consider  the subgroup of $ \aut_{\calge}(L)$ given by
 $$
 \aut_{\calge}^M(L)=\{\varphi\in\aut^M(L),\,[\varphi]\in \calge\},
 $$
 and denote by $
\calge^M\subset \cale^M(L)$
 the subgroup of relative homotopy classes of relative automorphisms whose free class live in $\calge$.

 In other terms, the subsequence
 $$
 \calge^M\to\calge\to\calge/H_0(L)\quad\text{of}\quad
\cale^M(L)\to \cale^*(L)\to \cale(L)
$$
 is identified with $G_\bq^{A_\bq}\to G_\bq^*\to G_\bq$, i.e.,  the sub sequence
 $$
\pi_0\aut_{G_\bq}^{A_\bq}(X_\bq)\to \pi_0  \aut_{G_\bq}^*(X_\bq)\to \pi_0  \aut_{G_\bq}(X_\bq)
\quad\text{of}\quad
 \cale^{A_\bq}(X_\bq)\to  \cale^*(X_\bq)\to  \cale(X_\bq).
 $$
 \end{definition}

\begin{definition}\label{explogrelativo}
Denote by $\derge L$ the connected sub dgl  of $\derr L$ defined by
$$\derge_kL=
 \left\{
 \begin{array}{cl}
 \der_kL & \text{ if } k>0,
 \\
\text{$\theta\in\ker D$  such that $e^{\theta}\in\aut_{\calge}(L)$} & \text{ if } k=0.
 \end{array}\right.
 $$
By \cite[Lemma 8.1]{fefuenmu0} this  dgl is of the type in (\ref{deresp}) and thus, see \cite[\S2]{fefuenmu1}, it is  complete. In particular $\derge_0L$ is complete. Equivalently, its  image by the exponential, which is precisely $\aut_{\calge}(L)$, is a Malcev complete group.

 In the relative case define the sub dgl $\derge^ML$ of $\der L$ by
 $$\derge_k^ML=
 \left\{
 \begin{array}{cl}
 \der_k^ML & \text{ if } k>0,
  \\
\text{$\theta\in\ker D$  such that $e^{\theta}\in\aut_{\calge}^M(L)$} & \text{ if } k=0.
 \end{array}\right.
 $$
Since  $\derge^M L\subset \derge L$ this dgl is also complete. As before, this implies that  $\aut_{\calge}^M(L)$ is a Malcev complete group.
 \end{definition}

Consider the twisted product
$$
L\timest\derge^ML
$$
in which both factors are sub dgl's and $[\theta,x]=\theta(x)$ for $x\in L$ and $\theta\in\derge^ML$. Note that, under the splitting (\ref{iso}), this is a subdgl of the  $\Hom(\quic(L),L)\timest\der L$ in (\ref{mainrelation}).
Then, we prove:
\begin{theorem}\label{main1} The realization of the cdgl fibration
$$
L\longrightarrow L\timest \derge^ML\longrightarrow\derge^ML
$$
 has the  homotopy type of the universal fibration
$$
X_\bq\longrightarrow B(*,\aut_{G_\bq}^{A_\bq}(X_\bq),X_\bq)\longrightarrow B\aut_{G_\bq}^{A_\bq}(X_\bq).
$$
 \end{theorem}
We first compute the homotopy fibre of the map $\langle L\rangle\to \langle L\timest \derge^ML\rangle$. To this end observe that the image of the map in Corollary \ref{seccion},
$$
\sigma\colon L\stackrel{\simeq}{\longrightarrow}(\Hom^{\overline\quic(M)}(\quic(L),L)\timest \der^ML,\widetilde D)
$$
is contained in $\Hom^{\overline\quic(M)}(\quic(L),L)\timest \derge^ML$ as $\ad_x^M\in \derge^ML$ for any $x\in L$. Consider then the commutative triangle
\begin{equation}\label{dia1}
\xymatrix{(\Hom^{\overline\quic(M)}(\quic(L),L)\timest \derge^ML,\widetilde D)\ar[rd]^\varrho& \\ L\ar[u]^\sigma\ar[r]&L\timest\derge^ML\\
}
\end{equation}
where, under the identification in (\ref{split})
$$
\Hom^{\overline\quic(M)}(\quic(L),L)\timest \der^ML\cong L\timest \Hom^M(\quic(L),L)\timest \der^ML,
$$
the map $\varrho$ is the surjective morphism which vanishes on $\Hom^M(\quic(L),L)$, is the identity on $\der^ML$ and $\varrho(x)=x+\ad_x^M$ for $x\in L$. Being a cdgl fibration, its realization,
\begin{equation}\label{restri}
\langle \varrho\rangle\colon \langle (\Hom^{\overline\quic(M)}(\quic(L),L)\timest \derge^ML,\widetilde D)\rangle\to \langle L\timest\derge^ML\rangle,
\end{equation}
is also a fibration with connected base and with fibre the realization of the  $\ker\varrho$,
\begin{equation}\label{restri2}
\langle (\Hom^M(\quic(L),L),D_q)\rangle.
\end{equation}
By Proposition \ref{modelofibrageneral} this is weakly equivalent to the subspace of $\map(X,X)$ consisting of maps extending the inclusion $i\colon A\hookrightarrow X$.
 Next, by restricting $\langle \varrho\rangle$ to the path component of the  trivial MC element we get a fibration  which fits in the commutative diagram,
$$
\xymatrix{\langle(\Hom^{\overline\quic(M)}(\quic(L),L)\timest \derge^ML,\widetilde D)\rangle^0\ar[rd]^(.56){\langle \varrho\rangle}& \\ \langle L\rangle \ar[u]^{\langle \sigma\rangle}_\simeq \ar[r]&\langle L\timest\derge^ML\rangle\\
}
$$
 resulting of realizing the path component of diagram  (\ref{dia1}) at 0. Thus, the homotopy fibre of $\langle L\rangle\to \langle L\timest \derge^ML\rangle$ is the fibre $F$ of $\langle \varrho\rangle$  in this diagram.

 \begin{lemma}\label{lemafibra}  $F={\amalg}_{[\gamma]\in\calge^M}\,\,\langle  \Hom^M(\quic(L),L)\rangle^{\gamma q}$.
 \end{lemma}

\begin{proof} Being $\langle\rho\rangle$ the restriction of (\ref{restri}), of fibre (\ref{restri2}), to a precise component, $F$ consists of all path components of $\langle (\Hom^M(\overline\quic(L),L,D_q)\rangle$ contained in the path connected complex $\langle(\Hom^{\quic(M)}(\quic(L),L)\timest \derge^ML,\widetilde D)\rangle^0$. Such a component, see Remark \ref{remim}, is determined by a class $[\gamma]\in\{L,L\}^M$ such that, $(\gamma-\id_L)q$ is gauge related to $0$ as MC elements in $(\Hom^{\overline\quic(M)}(\quic(L),L)\timest \derge^ML,\widetilde D)$. If this is the case via an element of degree zero in $\Hom^M(\overline\quic(L),L)$, it means that $\gamma\sim_M\id_L$.

On the other hand, if $\theta\in \derge^M_0L$ is such that $\theta\gauge(\gamma-\id_L)=0$, this translates to:
$$
e^{\ad_\theta}\bigl((\gamma-\id_L)q\bigr)-\frac{e^{\ad_\theta}-1}{\ad_\theta}(D\theta)=0.
$$
As $\theta$ is a cycle, this becomes
$
e^{\ad_\theta}\bigl((\gamma-\id_L)q\bigr)=0$ which reduces to $e^\theta\gamma=e^\theta$. Thus,
$
\gamma=e^{-\theta}\in \aut^\calge_M(L)$, that is, $[\gamma]\in \calge^M$.
\end{proof}
We then have a fibration sequence
$$
F\stackrel{\xi}{\longrightarrow}\langle L\rangle\longrightarrow \langle L\timest \derge^ML\rangle
$$
for which:
\begin{corollary}\label{cornue2}
$F\simeq \aut_{G_\bq}^{A_\bq}(X_\bq)$ and $\xi$ is null-homotopic.
\end{corollary}

\begin{proof}
By Theorem \ref{autofibra},
$$
\aut^{A_\bq}(X_\bq)\simeq \underset{[\gamma]\in\cale^M(L)}{\amalg} \langle \Hom^M(\quic(L),L)^{\gamma q}\rangle.
$$
By restricting this equivalence to the corresponding components, we have
$$
\aut_{G_\bq}^{A_\bq}(X_\bq)\simeq \underset{[\gamma]\in\calge^M}{\amalg} \langle \Hom^M(\quic(L),L)^{\gamma q}\rangle
$$
and the first assertion follows from Lemma \ref{lemafibra}.   For the second observe that  the map $\xi$ at a given component $\langle  \Hom^M(\quic(L),L)\rangle^{\gamma q}$ of $F$ is the realization of the path component at 0 of the composition
$$
(\Hom^M(\quic(L),L),D_{\gamma q})\to (\Hom^{\overline\quic(M)}(\quic(L),L)\timest \der^ML,\widetilde D)\stackrel{\eta}{\longrightarrow} L
$$
which is trivial.
\end{proof}

\begin{proof}[Proof of Theorem \ref{main1}]
Observe that the cdgl inclusion $M\timest \derge^ML\hookrightarrow L\timest \derge^ML$ is an $M$-section so that the sequence $
L\to L\timest \derge^ML\to\derge^ML
$ is an $M$-fibration.
As a consequence, see \S1.2, the sequence
$$
\langle L\rangle\longrightarrow \langle L\timest \derge^ML\rangle\longrightarrow\langle\derge^ML\rangle
$$
is an $A_\bq$-fibration. Moreover, by  Remark 8.5 of \cite{fefuenmu0}, the holonomy action of this fibration sequence $\pi_1\langle \derge^ML\rangle \to\cale\langle L\rangle$ is identified with the group morphism
$$
H_0\derge^ML\to \cale^M(L),\quad [\theta]\mapsto [e^\theta].
$$
However, $\theta\in \derge^M_0L$ if and only if  $e^{\theta}\in\aut_\calge^M(L)$ so that its relative homotopy class lies in $\calge^M$. Equivalently, the image of the holonomy action lies $G_\bq^{A_\bq}\subset \cale^{A_\bq}(X_\bq)$.
 Thus, see \S1.2,
by the classifying property of the corresponding universal fibration we have a commutative diagram
$$
\xymatrix{X_\bq\ar[r]&B(*,\aut_{G_\bq}^{A_\bq}(X_\bq),X_\bq)\ar[r]&B\aut_{G_\bq}^{A_\bq}(X_\bq)\\
\langle L\rangle\ar[u]^\simeq\ar[r]&
\langle L\timest\derge^M L \rangle\ar[u]\ar[r]&
\langle \derge^ML\rangle.\ar[u]\\}
$$
Then, by Corollary \ref{cornue2}, and taking into account that the restriction of the evaluation fibration $\ev\colon \aut(X_\bq)\to X_\bq$ to $\aut_{G_\bq}^{A_\bq}(X_\bq)$ is trivial, this diagram can be extended to the left as,
$$
\xymatrix{\aut_{G_\bq}^{A_\bq}(X_\bq)\ar[r]^(.56)\ev&
X_\bq\ar[r]&B(*,\aut_{G_\bq}^{A_\bq}(X_\bq),X_\bq)\ar[r]&B\aut_{G_\bq}^{A_\bq}(X_\bq)\\
F\ar[u]^\simeq\ar[r]_\xi&\langle L\rangle\ar[u]^\simeq\ar[r]&
\langle L\timest\derge^M L \rangle\ar[u]\ar[r]&
\langle \derge^ML\rangle.\ar[u]\\}
$$
As both rows are fibration sequences the last two vertical maps are equivalences.
\end{proof}

As a consequence we get:

\begin{theorem}\label{teofin} There is a homotopy commutative diagram
$$\begin{tikzcd}
\langle \derge^M L\rangle \arrow{d} \arrow{r}{\simeq}
&
B\aut^{A_\bq}_{G_\bq}(X_\bq) \arrow{d}
\\
\langle \derge L \rangle \arrow{r}{\simeq}
&
B\aut^{*}_{G_\bq}(X_\bq)
\end{tikzcd}$$
where the horizontal arrows are homotopy equivalences and both vertical arrows are induced by the obvious inclusions.
\end{theorem}

\begin{proof} Consider the following  diagrams
$$
\xymatrix{\langle L\timest\derge^M L \rangle\ar[r]\ar[d]_\simeq& \langle \derge^M L \rangle\ar[d]^\simeq&\langle L\timest\derge^M L \rangle\ar[r]\ar[d]& \langle \derge^M L \rangle\ar[d]\\
B(\ast,\aut^{A_\bq}_{G_\bq}(X_\bq),X_\bq)\ar[r]\ar[d]&B\aut^{A_\bq}_{G_\bq}(X_\bq)\ar[d]&\langle L\timest\derge  L \rangle\ar[r]\ar[d]_\simeq& \langle \derge  L \rangle\ar[d]^\simeq\\
B(\ast,\aut^{\ast}_{G_\bq}(X_\bq),X_\bq)\ar[r]&B\aut^{*}_{G_\bq}(X_\bq),&
B(\ast,\aut^{\ast}_{G_\bq}(X_\bq),X_\bq)\ar[r]&B\aut^{*}_{G_\bq}(X_\bq),\\}
$$
where the left top and bottom squares  are homotopy pullbacks by Theorems \ref{main1} and \cite[Proposition 7.8]{may} respectively. The right top square  is also a homotopy pullbacks since the realization functor preserves limits and fibrations. Finally, the right bottom square is a homotopy pullback by \cite[Theorem 0.2]{fefuenmu0}.
This exhibits the same map $\langle L\timest \derge^M L\rangle \to \langle \derge^M L\rangle$ as  homotopy pullbacks over the classifying fibrations
$$
X_\bq\to B(*,\aut_{G_\bq}^{*}(X_\bq),X_\bq)\to B\aut_{G_\bq}^{*}(X_\bq).$$

However, these homotopy pullbacks  are classified by homotopy classes of maps over $ B\aut_{G_\bq}^{*}(X_\bq)$ and thus, the composition of the right and left vertical arrows of both diagrams have to be homotopic.
\end{proof}
Consider the twisted product
$\derge L\timest sL
$
where $sL$ is an abelian sub Lie algebra and
$$
Dsx=-sdx+\ad_x,\quad[\theta,sx]=(-1)^{|\theta|}s\theta(x),\quad sx\in sL,\quad \theta\in\derge L.
$$
\begin{corollary}\label{corofin} The realization of the inclusion on the first factor
$$
\derge^ML\hookrightarrow \derge L\timest sL
$$
has the homotopy type of
$$
B\aut^{A_\bq}_{G_\bq}(X_\bq)\longrightarrow B\aut_{G_\bq}(X_\bq).
$$
\end{corollary}
\begin{proof}
Write the cdgl inclusion of the statement as the composition
$$
\derge^ML\hookrightarrow \derge L \hookrightarrow \derge L\timest sL
$$
whose realization, in view of Theorem \ref{teofin} and \cite[Theorem 0.1]{fefuenmu0}, has the homotopy type of the composition
$$
B\aut^{A_\bq}_{G_\bq}(X_\bq)\longrightarrow B\aut_{G_\bq}^*(X_\bq)\longrightarrow B\aut_{G_\bq}(X_\bq),
$$
which is $
B\aut^{A_\bq}_{G_\bq}(X_\bq)\longrightarrow B\aut_{G_\bq}(X_\bq)$.
\end{proof}

Consider the following dgl sequence formed by inclusions and projections,
\begin{equation}\label{muysimple}
\begin{aligned}
\Hom^M&(\quic(L),L)\to  \Hom(\quic(L),L)\to\Hom(\quic(M),L)\to\\
&\to\Hom(\quic(M),L)\timest\derge L\timest sL\to \derge L\timest sL.
\end{aligned}
\end{equation}
We finish by exhibiting the  homotopy type of the fibration sequence in Theorem \ref{prin3} as the realization of the appropriate components of the above sequence:

\begin{theorem}\label{final1}
 The sequence of fibrations
$$
\aut^{A_\bq}_{G_\bq}(X_\bq)\longrightarrow \aut_{G_\bq}(X_\bq)\stackrel{i_\bq^*}{\longrightarrow}\fmap_{G_\bq}(A_\bq,X_\bq)\longrightarrow B\aut^{A_\bq}_{G_\bq}(X_\bq)\longrightarrow B\aut_{G_\bq}(X_\bq)
$$
has the homotopy type of
$$
\begin{aligned}
\underset{[\gamma]\in\calge^M}{\amalg} \langle \Hom^M&(\quic(L),L)^{\gamma q}\rangle \to\underset{[\varphi]\in\calge/H_0(L)}{\amalg}\langle \Hom(\quic(L),L)^{\varphi q}\rangle \to\underset{[\varphi]\in\calge/H_0(L)}{\amalg}\langle \Hom(\quic(M),L)^{\varphi q}\rangle\\
&\to\underset{[\varphi]\in\calge/H_0(L)}{\amalg}\langle (\Hom(\quic(M),L)^{\varphi q}\timest\derge L\timest sL)^{\varphi q}\rangle\to\langle\derge L\timest sL\rangle.
\end{aligned}
$$
\end{theorem}

\begin{proof} The restriction of Theorem \ref{autofibra} to the corresponding components proves the assertion for the first fibration in the sequence. It remains to show that
$$
\fmap_{G_\bq}(A_\bq,X_\bq)\longrightarrow B\aut^{A_\bq}_{G_\bq}(X_\bq)\longrightarrow \aut_{G_\bq}(X_\bq)
$$
has the homotopy type of
$$
\begin{aligned}
&\underset{[\varphi]\in\calge/H_0(L)}{\amalg}\langle \Hom&(\quic(M),L)^{\varphi q}\rangle\to\underset{[\varphi]\in\calge/H_0(L)}{\amalg}\langle (\Hom(\quic(M),L)^{\varphi q}\timest\derge L\timest sL)^{\varphi q}\rangle\to\\&\to\langle\derge L\timest sL\rangle.
\end{aligned}
$$
To this end factor the inclusion $
\derge^ML\hookrightarrow \derge L\timest sL$, which models $B\aut^{A_\bq}_{G_\bq}(X_\bq)\to \aut_{G_\bq}(X_\bq)$, as
$$
\xymatrix{(\Hom(\quic(M),L)\timest \derge L\timest sL,D_{jq})\ar[rd] \\ \derge^ML\ar[u]^\kappa_\simeq\ar[r]&\derge L\timest sL\\
}
$$
where $\kappa$ is, as in Proposition \ref{proponueva}, the inclusion in the second factor.
As the diagonal map is a cdgl fibration, its realization
$$
\langle (\Hom(\quic(M),L)\timest \derge L\timest sL,D_{jq})\rangle\to \langle \derge L\timest sL\rangle
$$
is also a fibration with fibre
$$
\langle (\Hom(\quic(M),L),D_{jq})\rangle.
$$
By restricting this fibration to the component of $0$ we get a commutative diagram
$$
\xymatrix{\langle (\Hom(\quic(M),L)\timest \derge L\timest sL,D_{jq})\rangle^0\ar[rd] \\ \langle \derge^ML\rangle \ar[u]^{\langle\kappa\rangle}_\simeq\ar[r]&\langle \derge L\timest sL\rangle\\
}
$$
exhibiting the homotopy fibre of $ B\aut^{A_\bq}_{G_\bq}(X_\bq)\to \aut_{G_\bq}(X_\bq)$ as the fibre of the diagonal fibration. Finally, an analogous argument to the one in the proof of Lemma \ref{lemafibra} shows that this fibre is precisely $\underset{[\varphi]\in\calge/H_0(L)}{\amalg}\langle \Hom(\quic(M),L)^{\varphi q}\rangle$.
\end{proof}

The pointed version of (\ref{muysimple}) is the sequence
$$
\Hom^M(\quic(L),L)\to  \Hom(\overline\quic(L),L)\to\Hom(\overline\quic(M),L)
\to\Hom(\overline\quic(M),L)\timest\derge L\to \derge L.
$$
 Invoking now the model of the pointed fibration in Theorem \ref{autofibra} and Proposition \ref{varsigmareducida}, the arguments in the proof of Theorem \ref{final1} can be followed mutatis mutandis to show its pointed version:

 \begin{theorem}\label{final2}
 The sequence of fibrations
$$
\aut^{A_\bq}_{G_\bq}(X_\bq)\longrightarrow \aut_{G_\bq}^*(X_\bq)\stackrel{i_\bq^*}{\longrightarrow}\fmap_{G_\bq}^*(A_\bq,X_\bq)\longrightarrow B\aut^{A_\bq}_{G_\bq}(X_\bq)\longrightarrow B\aut_{G_\bq}^*(X_\bq)
$$
has the homotopy type of
$$
\begin{aligned}
\underset{[\gamma]\in\calge^M}{\amalg} \langle \Hom^M&(\quic(L),L)^{\gamma q}\rangle \to\underset{[\varphi]\in\calge(L)}{\amalg}\langle \Hom(\overline\quic(L),L)^{\varphi q}\rangle \to\underset{[\varphi]\in\calge(L)}{\amalg}\langle \Hom(\overline\quic(M),L)^{\varphi q}\rangle\\
&\to\underset{[\varphi]\in\calge(L)}{\amalg}\langle (\Hom(\overline\quic(M),L)^{\varphi q}\timest\derge L)^{\varphi q}\rangle\to\langle \derge L\rangle.
\end{aligned}
$$
\hfill$\square$
\end{theorem}

\section{Some examples}

We begin with an elementary one.

\begin{ex}[Boundary of a disk]

Consider the inclusion $A=S^n\hookrightarrow D^{n+1}=X$ of the $n$-sphere into the $n+1$-disk, with $n\ge 1$. The Lie minimal model of this cofibration is $M=(\mathbb{L}(x),0)\hookrightarrow L=(\libc(x,y),d)$, where $|x| =n-1, |y|=n$ and $dy=x$. Here, both  $\aut^{A}(X)$ and $\aut^{A_\bq}(X_\bq)$ are contractible and,  in particular, $\cale^A(X)=\cale^{A_\bq}(X_\bq)$ reduces to the identity. This amounts to say that any cdgl automorphism $\varphi\colon L\to L$ with $\varphi(x)=x$ and $\varphi(y)=y+\Psi$, with $\Psi\in L^2$, has to be homotopic to the identity. In the case of $n>1$ this is trivial as $\varphi(y)=y$. However for $n=1$ we have an interesting phenomenon: in this case
$$\varphi(y)=y+\sum_{i\geq 1} \lambda_i\ad_x^i(y),\quad \lambda_i\in\bq.$$
Finding an explicit homotopy of the form $\phi\colon L\to \wedge(t,dt)\widehat \otimes L$ between $\varphi$ and the identity could be cumbersome. Nevertheless, we can translate this problem to the world of derivations and think of this process as a ``linearization'' of the problem. As $B\aut^{A_\bq}(X_\bq)$ is contractible, $\derge^ML$ is acyclic in view of Theorem \ref{main1}. In particular,  any  derivation $\theta\in \derge_0^ML$, which is necessarily   of the form
$$\theta(x)=0, \;\;\; \theta(y)=\sum_{i\geq 1} a_i \ad_x^i(y),\quad a_i\in\bq,$$
is the boundary of some $\eta\in\derge_1^M(L)$. Write
$$\eta(x)=0, \;\;\; \eta(y)=\sum_{i\geq 0} b_i [y,\ad_x^i(y)]$$
so that a short computation provides
$$(D\eta)(x)=0, \;\;\; (D\eta)(y)= 2b_0 [x,y]+ \sum_{i\geq 1}b_i \ad_x^{i+1}(y).$$
We solve the equation by setting $b_0=a_1/2$ and  $b_i=a_{i+1}$, for $i\ge 1$.
\end{ex}

In the remaining of this section we present a general framework covering a wide spectrum of examples. The following complements  the work in  \cite{bens}.
\medskip

Attach an $(n+1)$-cell, $n> 0$, to the  finite nilpotent CW-complex $A$ through the map $f\colon S^{n}\to A$ so that the resulting complex $X=A\cup_fe^{n+1}$ is nilpotent. If $M=(\libc(U),d)$  is a Lie model of $A$, a Lie model of $X$ is the extension of $M$ given by
$$
L=(\libc(U\oplus\bq x),d),\quad |x|=n, \quad dx=\omega,
$$
with $[\omega]\in H_*(\libc(U),d)$ representing $[f_\bq]\in\pi_n(A_\bq)$. Note that if $f_\bq$ is homotopically trivial we may, and will assume $x$ to be a cycle.

\begin{proposition}\label{propofin1} For any  $k\ge 2$,
$$\pi_k B\aut^{A_\bq}(X_\bq)\cong \pi_{n+k}(X_\bq).
$$
\end{proposition}

\begin{proof}
An easy inspection shows  that the map
\begin{equation}\label{eq:bijection}
\psi\colon \Der_*^ML\stackrel{\cong}{\longrightarrow} L_{*+n},\;\;\; \theta\mapsto \theta(x),
\end{equation}
is an isomorphism of chain complexes. The result follows by recalling that $\derr_{k\ge1}L=\derge_{k\ge 1} L $ for any given $G$.
\end{proof}
Note that this covers in particular Example 5.1 of \cite{bersa} where the considered case is the inclusion  $\C P^k\hookrightarrow \C P^n$, with $1\le k<n$.

\medskip

 In the next result, we also describe the group
$$\pi_1B\aut^{A_\bq}(X_\Q)\cong  \cale^{A_\bq}(X_\bq).
$$

\begin{rem}\label{finalrem}
Observe that the isomorphism in (\ref{eq:bijection}) induces, in general, a new dgl structure in $L$ with a different grading. Moreover, as we will see, $H_n(L)$ if $f_\bq$ is not homotopically trivial, or a codimension $1$ subspace of  $H_n(L)$ otherwise,  inherit  structures of nilpotent rational groups from $H_0(\derr L)$ with the BCH product.
\end{rem}

\begin{theorem}\label{teoejemplo}
There is a nilpotent group structure on $\pi_{n+1}(X_\bq)$ for which
$$
\cale^{A_\bq}(X_\bq)\cong \pi_{n+1}(X_\bq)
$$
if $f_\bq$ is not homotopically trivial.
Otherwise, there is a short exact sequence of groups
$$
K\to \cale^{A_\bq}(X_\bq)\to\bq^*
$$
where $K$ is a  a codimension $1$ subspace of $ \pi_{n+1}(X_\bq)$ endowed with a rational nilpotent group structure.
\end{theorem}

\begin{proof}
Let $\Gamma\subset H_*(X)$ be the image of the map $H_*(A)\to H_*(X)$ induced by the inclusion. This subspace defines a filtration
$$ H_*(X)\supset \Gamma \supset \{0\}$$
and we consider ${G}\subset \cale(X)$  the subgroup which stabilizes the filtration above, i.e. $[f]\in \cale(X)$ belongs to ${G}$ if and only if
$$
(H_*(f)-\id)(H_*(X))\subset \Gamma  \text{ and } (H_*(f)-\id)(\Gamma)=0.$$
By definition this group acts nilpotently on $H_*(X)$.
Further, $G^A\subset \cale^A(X)$  consists of relative homotopy classes of self equivalences $f$ of $X$ which fix $A$ and such that $(H_*(f)-\id)(H_*(X))\subset \Gamma$, i.e.,   $H_*(f)$ is the identity on $H_*(X)/\Gamma$. By \cite[Theorem 7.6]{fefuenmu0} $G_\bq$ is the subgroup of $\cale(X_\bq)$ which stabilizes the filtration
$$ H_*(X_\bq)\supset \Gamma_\bq \supset \{0\}$$ so that again, $G_\bq^{A_\bq}\subset \cale^{A_\bq}(X_\bq)$ consist of those relative homotopy classes of sel equivalences inducing the identity on $H(X_\bq)/\Gamma_\bq$.

At the Lie side, $H_*(A_\bq)\to H_*(X_\bq)$ is identified with the map
$$
U\longrightarrow H_*(U\oplus \bq x,d_1)
$$
where $d_1U=0$ and $d_1x$ is the linear part of $\omega$. This is obviously the inclusion $U\to U\oplus \bq x$ if $d_1x=0$ or the projection $U\to U/(d_1x)$ otherwise.

Hence,
\begin{equation}\label{aver}
\calge^M=\{[\varphi]\in \cale^M(L) \text{ such that } \varphi_1(x)-x\in U\},
\end{equation}
and thus, one easily sees that
$$\derge_0^ML=\{\theta\in \Der_0^ML\mid D\theta=0, \; \theta_1(x)\in U \}.$$
Therefore, using (\ref{eq:bijection}), $\derge_0^ML$ corresponds to the subspace $L'_n=\ker d\cap L''_n\subset L_n$ where $L''_n$ is the complement of the subspace $\bq x$. An easy inspection shows that $\ker d\subset L''_n$ if and only if $dx$ is not a boundary on $M$, that is, if the ataching map $f_\bq$ is not homotopically trivial. If this is the case
$$
\dim H_0(\derge^ML)=\dim H_n(L)-1.
$$
 Otherwise,
 $$H_0(\derge^ML)\cong H_n(L).
 $$

On the other hand, and element $\varphi\in \aut^M(L)$ is completely determined by its image on  $x$. A short computation show that, $\varphi(x)=\lambda x+\alpha$ with $\lambda\in \bq^*$ and $\alpha$ a cycle in $L''_n$. Moreover, if $dx$ is not a boundary, then $\lambda=1$. In view of (\ref{aver}), this amounts to say that
$$\cale^M(L)=\calge^M.
$$

Conversely, if the attachment is trivial, any $\lambda\neq 0$ and any cycle $\alpha$ produces an isomorphism $\varphi$ and we have  a group inclusion
$$\aut_{\mathcal{G}}^M(L)\to \aut^M(L)$$
whose cokernel is $\bq^*$. Since $\aut_1^M(L)\subset \aut_{\mathcal{G}}^M(L)$, taking quotients in the above inclusion provides a  short exact sequence of groups
$$\calge^M\to \cale^M(L)\to \bq^*.$$

By all of the above, in view of  Theorem \ref{main1} and Remark \ref{finalrem}, and through the isomorphisms
$$\cale^M(L)\cong \cale^{A_\bq}(X_\bq),\quad \calge^M\cong G_\bq^{A_\bq}\cong \pi_1B\aut^{A_\bq}_{G_\bq}(X_\Q)\cong H_0(\derge^ML),
$$
the theorem follows choosing $K=\cale^{A_\bq}_{G_\bq}(X_\bq)$.

\end{proof}

\eject

\noindent {\sc Institut de Math\'ematiques et Physique, Universit\'e Catholique de Louvain, Chemin du Cyclotron 2,
1348 Louvain-la-Neuve,
         Belgique}.

\noindent\texttt{yves.felix@uclouvain.be}

\bigskip

\noindent{\sc Centro de Investigaci\'on en Matem\'aticas, Unidad de M\'erida, Sierra Papacal, M\'erida, YUC 97302 M\'exico.}

\noindent
\texttt{mario.fuentes@cimat.mx}

\bigskip

\noindent{\sc Departamento de \'Algebra, Geometr\'{\i}a y Topolog\'{\i}a, Universidad de M\'alaga, 29080 M\'alaga, Spain.}

\noindent
\texttt{aniceto@uma.es}

 \end{document}